\documentclass[12pt]{amsart}
\usepackage{amsthm,amsmath,amssymb}
\usepackage{addfont}
\usepackage{mathrsfs}
\usepackage[T1]{fontenc}
\newcommand{\n}{\noindent}
\numberwithin{equation}{section}

\def\ca{{\mathcal A}}

\def\cq{{\mathcal Q}}
\def\car{{\mathcal R}}
\def\cas{{\mathcal S}}
\def\ct{{\mathcal T}}

\def\cv{{\mathcal V}}

\def\bc{{\mathbb C}}

\def\bn{{\mathbb N}}

\def\br{{\mathbb R}}

\def\a{\alpha}

\def\b{\beta}
\def\d{\delta}

\def\g{\gamma}
\def\k{\kappa}
\def\l{\lambda}       
\def\m{\mu}
\def\n{\nu}

\def\f{\varphi}

\theoremstyle{plain}
\newtheorem{lemma}{Lemma}[section]
\newtheorem{proposition}[lemma]{Proposition}
\newtheorem{theorem}[lemma]{Theorem}
\newtheorem{corollary}[lemma]{Corollary}
\theoremstyle{definition}

\newtheorem{definition}[lemma]{Definition}

\begin{document}

\title[Grothendieck's inequalities  ]{\textsc{  Some points of view on Grothendieck's inequalities}}

\author[E.~Christensen]{Erik Christensen}
\address{\hskip-\parindent
Erik Christensen, Mathematics Institute, University of Copenhagen, Copenhagen, Denmark.}
\email{echris@math.ku.dk}
\date{\today}
\subjclass[2010]{ Primary: 15A39, 46B28, 46L07. Secondary: 15A60, 46B85, 47L25.}
\keywords{ Grothendieck inequality, bilinear operators,  Schur product,  completely bounded,  operator space, Stinespring representation, duality, tensor product  } 

\begin{abstract}   
 Haagerup's proof of the {\em non commutative little Groth-endieck inequality}  raises some questions on the {\em commutative little inequality,} and it offers a new result on scalar matrices with non negative entries. The theory of completely bounded maps may be used to show that the {\em commutative Grothendieck inequality }  follows from the {\em little commutative  inequality,}  and that this passage  may be given a geometric form as a relation between a pair of compact convex sets of positive matrices, which, in turn, characterizes the {\em little constant $k_G^\bc.$}        
 \end{abstract}

\maketitle

\section{Introduction}

Grothendieck's work on tensor products of normed  spaces \cite{Gr} has influenced mathematics in several ways, some of which are very surprising. This is described in Pisier's survey article \cite{Pi3} and the book \cite{DHS} by Diestel, Fourie and Swart.  Here we will focus on the inequalities which are named {\em Grotthendiek's inequality} and { \em Grothendieck's little inequality} in the setting of complex $m \times n$ matrices. This is a  continuation  of our recent articles \cite{C3}, \cite{C4}, where we showed that the theory of operator spaces and completely bounded maps provides a set up, which fits very well - in our opinion - to the existing results related to Grothendieck's inequalities.   An important aspect in Grothendieck's work deals with a bounded operator between Banach spaces which factors through a Hilbert space or through the space of continuous complex functions on a compact  topological space. In this article the Hilbert spaces are finite dimensional and the compact spaces we will meet have only finitely many points, so our results will deal with complex $m \times n$ matrices, the set of which we denote $M_{(m,n)}(\bc).$ 
This set of matrices is in a canonical way isomorphic to the algebraic tensor product $\bc^m \otimes \bc^n $ where the isomorphism  is described via the canonical basis $(\d_i)_{\{1 \leq i \leq m\}}$ for  $\bc^m$ and $(\g_j)_{\{1 \leq j\leq n\}},$ for $\bc^n$ and the matrix units $\{e_{(i,j)} : 1 \leq i \leq m,\, \, 1 \leq j \leq n\}$ for $M_{(m,n)}(\bc)$ by  the linear  map $\f : \bc^m \otimes \bc^n \to M_{(m,n)}(\bc ) $ which satisfies $\f(\d_i \otimes \g_j ):= e_{(i,j)}.$
In several spots we will use that the image $\f(a \otimes b) $ is the rank one matrix with entries $\f(a \otimes b)_{(i,j)}= a_ib_j$ and also use, that this matrix is a product of a one column matrix $a_| := (a_1, \dots, a_m)$ and a one row matrix $b_{-} := (b_1, \dots, b_n),$ so $\f(a \otimes b ) = a_|b_{-}.$
 
 For each real positive $p \geq 1 $ or $p = \infty$ and any natural number $k$ we let $\|.\|_p$ denote the usual $p-$norm on $\bc^k.$ Given a couple of normed spaces   such as $(\bc^m, \|.\|_p)$ and $(\bc^n, \|.\|_r), $ we recall that Schatten \cite{Scha} has introduced the concept named a cross norm on the tensor product $(\bc^m, \|.\|_p) \otimes (\bc^n, \|.\|_r),$ and we recall that a norm say $|||.|||$  on the tensor product of the normed spaces is called a cross norm on this tensor product of normed spaces if it satisfies $$ \forall \eta \in \bc^m \, \forall\xi \in \bc^n : \quad |||(\eta \otimes \xi)|||  = \|\eta\|_p \|\xi\|_r.$$ 
Schatten proved that there is a minimal and a maximal cross norm. Today the minimal cross norm  is called the {\em injective cross norm } and denoted $\|.\|_\vee. $ The maximal cross norm is called the {\em projective cross norm } and it is denoted $\|.\|_\wedge.$ In the situation with $(\bc^m, \|.\|_p) \otimes ( \bc^n, \|.\|_r)$  we can then define  norms $\|.\|_{\vee(p,r)}$ and $\|.\|_{\wedge(p,r)} $ on $M_{(m,n)}(\bc) $ by a transport of the injective and projective norms on the tensor product to norms on $M_{(m,n)}(\bc)$ via the isomorphism $\f$ we defined above. 
 
 There are many well known norms on $M_{(m,n)}(\bc),$ and amongst them we will right now  mention the {\em operator norm, } which we denote $\|X\|_\infty = \|X|_{\vee(2,2)},$ the {\em Hilbert Schmidt norm,}  which we denote $\|X\|_2$ and the {\em Schur multiplier norm,} which we denote $\|X\|_S.$ We will remind you on the Schur product of matrices. Given two complex matrices $X = (X_{(i,j)}) $ and $A= (A_{(i,j)} )$ in $M_{(m,n)} (\bc),$ then we define their Schur product $X \circ A$ to be the matrix in $M_{(m,n)}(\bc) $ which is given by the equation $(X \circ A)_{(i,j)} := X_{(i,j)}A_{(i,j)}.$ 
 We can then formulate Grothendieck's inequality in a way which is close to the original one, except - of course - for the use of Grothendieck's name.   
 \begin{theorem} \label{BigK}
There exists a positive real $K_G^\bc \leq \sinh(\pi/2)$ such that for any complex  $m \times n $ matrix $X$ we have $\|X\|_{\wedge(\infty, \infty)} \leq K_G^\bc\|X\|_S.$ 
\end{theorem}  
The exact value of $K_G^\bc$ is unknown, but after the combined efforts of several authors Pisier reports in Section 4 of  \cite{Pi3} that $1.338 < K_G^\bc \leq 1.4049.$ The Grothendieck inequality is most often described as a property for bilinear forms  on the product $(\bc^m, \|.\|_\infty )\times (\bc^n, \|.\|_\infty)$ and we will return to this formulation of the inequality later on.

 There is also an inequality which often is named {\em Grothendieck's little inequality.} To  formulate that one we recall the norms $\|X\|_F$ and $\|X\|_{cbF} $  from   \cite{C3},  which we defined on $M_{(m,n)}(\bc).$ The norm  $\|X\|_F$  is defined as the norm of the linear operator $F_X$ with the matrix $X$ acting as an operator from $(\bc^n, \|.\|_\infty)$ to $(\bc^m, \|.\|_2),$ and then $\|X\|_{cbF}$ is the completely bounded norm of $F_X.$ It follows from the definition of the injective norm that $\|X\|_F = \|X\|_{\vee(2,1)}.$ The little inequality may then be formulated as follows.

 \begin{theorem} \label{Litk} 
 There exists a positive constant $k_G^\bc$ such that for any complex $m \times n$ matrix $X$ we have   $$\|X\|_{cbF}  \leq \sqrt{k_G^\bc} \|X\|_F = \sqrt{k_G^\bc} \|X\|_{\vee(2,1)}.$$
 \end{theorem}

 It is known, see Section 5 of  \cite{Pi3}, that $k_G^\bc = 4/\pi.$ 
 
  The normed space  $(\bc^n, \|.\|_\infty ) $  may be considered to be the $n$ dimensional abelian C*-algebra $\ca_n := C(\{1, \dots, n\}, \bc),$ the continuous complex functions on the set $\{1, \dots, n\}$ equipped with the sup-norm, and in this way Grothendieck's inequality contains  a statement on bounded  bilinear forms on a product of 2 abelian C*-algebras. This raises the natural question if Grothendieck's inequality does have an extension to bounded bilinear forms on a pair of non commutative C*-algebras.
   This was solved by Pisier in \cite{Pi4}, where he shows what the non commutative inequality ought to be and also shows that this inequality holds if a certain approximation property holds. This approximation   restriction was not a serious  problem and it was removed by Haagerup in \cite{Ha2}. 
   We will not go into any discussion of the content of the non commutative Grothendieck inequality here, but mention that Haagerup in an appendix to \cite{Ha1}  gives a proof of the non commutative little Grothendieck inequality which actually seems to contain new information when applied in the finite dimensional and abelian situation we are studying here. This aspect is discussed in Section \ref{Ha} below. Another aspect of this application of Haagerup's work is, that we can prove Grothendieck's little inequality from Theorem \ref{Litk} with elementary mathematics, but at a cost of the  rather bad inequality $k_G^\bc \leq 2.$ 
   
   Our research in the articles \cite{C3} and \cite{C4} studied relations between some norms on $M_{(m,n)}(\bc)$ and it gave some characterizations of the norms in terms of certain factorization properties. 
   The norms we studied are not new, but the perspective is to look at them as   completely bounded norms of some linear or bilinear operators and then investigate their minimal Stinespring representations, as defined in Definition 2.1 of \cite{C2}.  The concepts named   {\em completely bounded} and {\em Stinespring representation} come from the theory of {\em operator spaces and completely bounded maps,}  but  we will not give an introduction to that theory here. We gave a short description of the most needed facts we use, in \cite{C3} on page 546, right after Proposition 1.5 of that article,  and there are fine text books by Paulsen \cite{Pa}, Pisier \cite{Pi2} and Effros and Ruan \cite{ER} which describe this subject.
   
    Already Grothendieck introduced many norms in his résumé \cite{Gr}, Pisier in \cite{Pi3} and  Diestel, Fourie and Swart in \cite{DHS} list more norms than we will discuss here.  We will look at  some of the norms mentioned in Section 3 of  Pisier's survey, but there is   also the norm $\|.\|_T$  defined in Definition 1.2 of \cite{C4}, which may exist under a different name in the literature ? 
   
   In Section 2 we recall the norms we studied in \cite{C3} and \cite{C4} as factorization norms  on matrices in $M_{(m,n)}(\bc).$ It is not obvious that they are all cross norms, we think, but the factorization result of \cite{C4} gives an easy way to verify this. The results from  Section 3 of Pisier's survey makes it easy to identify all but one  of the norms we have introduced  with norms from \cite{Pi3}.

In Section 3 we go back to Haagerup's proof in \cite{Ha1}  of the little Grothendieck inequality for non commutative C*-algebras, and we find, that when his construction is applied in the finite dimensional and abelian setting we are studying, then the objects,  he investigates, become  quite concrete and in principle computable. This raises the question: Will this method of Haagerup's actually  produce the optimal cbF-factorization or operator factorization of $X$ which we studied in item (i) of Theorem 2.1 of \cite{C4}. We guess that the answer in general  is no, but we prove that the answer is yes for matrices with non negative entries only. For such matrices we actually show that Grothendieck's little inequality holds with the constant equal to 1.  
  
  In Section 4 we return to \cite{C4}, where we showed that for a positive matrix $P$ in $M_n(\bc) $ we have $\|P\|_{cbB} \leq k_G^\bc\|P\|_B,$ or in words that for positive matrices the constant $K_G^\bc$ may be replaced by $k_G^\bc.$   When applying this result to a certain  positive matrix in $M_{(m+n)}(\bc) $ we can obtain  bounds for the constant $K_G^\bc $  - based on the value of $k_G^\bc$ - as \begin{equation} \label{bounds} k_G^\bc \leq K_G^\bc \leq  \ k_G^\bc/(2 - k_G^\bc).\end{equation}
  
  The Theorem 1.1 is  dual to the Grothendieck inequality and it is really a statement on the convex hull of the rank one matrices $\f(\eta\otimes \xi) $ with  $\eta \in \bc^m,$ $\|\eta\|_\infty \leq 1,$ $\xi \in \bc^n$ and  $\|\xi\|_\infty \leq 1. $ This may be seen as a geometrical formulation of Grothendieck's inequality. We pursue this aspect, but for positive matrices, and we show in Theorem \ref{geo} that the mentioned result  for positive matrices "$\|P\|_{cbB} \leq k_G^\bc\|P\|_B$" implies a certain property for the closed convex hull of the positive rank one matrices with diagonal equal to the identity. That geometrical result turns out to characterize $k_G^\bc$ and it has  equation \eqref{bounds} as  an easy corollary.

\section{Cross norms and different names} \label{cross}
Our identification of $\bc^m \otimes \bc^n $ with $M_{(m,n)}(\bc)$ goes via the linear map $\f$ given by $\f(\eta \otimes \xi ) = \sum_{i,j} \eta_i\xi_j e_{(i,j)} .$ From here you see that the image by $\f$ of the non vanishing   elementary tensors $\eta \otimes \xi,$ consists of the rank one matrices.    Hence in order to establish that a norm on $(\bc^m, \|.\|_p) \otimes (\bc^n, \|.\|_r)$ is a cross norm we will have to look at rank 1 matrices only. There are quite a few norms 
involved  in this section. From Section 3 of \cite{Pi3} we get some projective and some injective norms, the norms $\|.\|_H,$ $\|.\|_{H'} ,$  the Schur multiplier norm and the Haagerup tensor norm, which we denote $\|.\|_h.$ From Definition 1.2 of \cite{C4} we  get  norms with subscripts $\{ F, cbF, B, cbB, S, T\, \}.$   
The Section 3 of \cite{C3} gives a characterization of the norms $\|.\|_{cbF},$ $  \|.\|_{cbB}, \|.\|_S$ and $ \|.\|_T$ in terms of certain factorization properties, and based on theses properties we can obtain concrete factorizations of rank one matrices, which show that all the norms are  cross norms on tensor products of the form $(\bc^m, \|.\|_p) \otimes (\bc^n, \|.\|_r)$ for some choices of $p$ and $r$ in the set $\{1,2, \infty\}.$  When possible, we  will  identify some of the norms Pisier mentions with some of the norms we have given different names. The reason why we have introduced the new names is, that we find that these names seem to fit well with our {\em completely bounded } approach to the problems under investigation.

\begin{theorem}
Let $X$ be a non-zero complex $m \times n$ matrix  of rank 1 with $x_{(i,j)} = \m_i\n_j$ for vectors $\m$ in  $\bc^m$ and $\n$ in $\bc^n,$ then \begin{itemize}
\item[(i)] $\|X\|_F = \|X\|_{cbF} = \|\m\|_2\|\n\|_1.$
\newline The norm $\|.\|_F$ is the minimal cross norm $\vee(2,1) $  on \newline  $(\bc^m, \|.\|_2) \otimes (\bc^n, \|.\|_1).$ \newline The norm $\|.\|_{cbF}$  is a cross norm on   on  $(\bc^m, \|.\|_2) \otimes (\bc^n, \|.\|_1),$ and it satisfies $\|.\|_{cbF} \leq \sqrt{k_G^\bc} \|.\|_F.$ The norm $\|.\|_{cbF} $ is conjugate dual to the norm $\|.\|_T.$ 

\item[(ii)] $\|X\|_B = \|X\|_{cbB} = \|\m\|_1\|\n\|_1.$ \newline The norm $\|.\|_B $ is the minimal  cross norm
$\vee(1,1)$  on  \newline $(\bc^m, \|.\|_1)\otimes (\bc^n, \|.\|_1).$ \newline The norm $\|.\|_{cbB}$  is a cross norm  on  $(\bc^m, \|.\|_1)\otimes (\bc^n, \|.\|_1),$ and it satisfies $\|. \|_{cbB} \leq K_G^\bc \|.\|_B.$The norm $\|.\|_{cbB} $ is conjugate dual to the norm $\|.\|_S.$ 

\item[(iii)] $\|X\|_S = \|\m\|_\infty\|\n\|_\infty.$
\newline The norm $\|.\|_S = \|.\|_{cbS}$ is a  cross norm  on \newline  $(\bc^m, \|.\|_\infty) \otimes (\bc^n, \|.\|_\infty).$ It satisfies $\|.\|_{\wedge(\infty, \infty)} \leq K_G^\bc\|.\|_S.$ The norm $\|.\|_{S} $ is conjugate dual to the norm $\|.\|_{cbB}.$  
\item[(iv)] $\|X\|_{T} = \|\m\|_2\|\n\|_\infty.$ \newline The norm $\|.\|_T = \|.\|_{cbT}$ is a cross norm on \newline $(\bc^m, \|.\|_2) \otimes (\bc^n, \|.\|_\infty). $ It satisfies $\|.\|_{\wedge(2, \infty) } \leq \sqrt{k_G^\bc} \|.\|_T.$ The norm $\|.\|_{T} $ is conjugate dual to the norm $\|.\|_{cbF}.$ 
\item[(v)] With the notation from Section 3 of \cite{Pi3}, $\|.\|_S = \|.\|_H = \|.\|_h$.
\item[(vi)] With the notation from Section 3 of \cite{Pi3}, $\|.\|_{cbB} = \g_2(.)= \|.\|_{H'} $. 
\end{itemize}
The different factorizations  of $X,$ which will show that  the cb-norms are as claimed, are described in the proof of the proposition.
\end{theorem}

\begin{proof}
It follows from a direct computation that $\|X\|_F = \|\m\|_2\|\n\|_1.$  We will show that the vector $\xi$ in the cbF-factorization of $X$  is given by $\xi_j :=( |\n_j|^{(1/2)})/(\|\n\|_1^{(1/2)}).$ By construction this $\xi$ is a positive unit vector in $(\bc^n, \|.\|_2).$ We may define $\Delta_n(\xi)^{inv}$ as the diagonal matrix with entries equal to $\xi_j^{-1} $ if $\xi_j >0$ and $0$ if $\xi_j = 0.$ Then  we have  $\|X\Delta_n(\xi)^{inv}\|_\infty = \|X\Delta_n(\xi)^{inv}\|_2$ since the rank is one.  Then this norm is given as $$\|X\Delta_n(\xi)^{inv}\|_\infty = \sqrt{ (\sum_i |\mu_i|^2 )(\sum_j |\n_j|\|\n\|_1)} = \|\m\|_2 \|\n\|_1.$$  The by Theorem 1.3 item (i) of \cite{C4}  we have $$\|X\|_{cbF} \leq \|X\Delta_n(\xi)^{inv}\|_\infty = \|\m \|_2\|\n\|_1 = \|X\|_F \leq \|X\|_{cbF},$$ so the norms are cross norms as claimed. The very definition of the norm $\|.\|_F$ tells that it is the minimal cross norm $\vee(2,1).$  The inequality $\|.\|_{cbF} \leq \sqrt{k_G^\bc} \|.\|_{\vee(2, 1)}$ then follows from Theorem 2.1 of \cite{C3}. The duality statement in item (i) is a consequence of equation (4.4) of \cite{C3}, and item (i) follows.   

 With respect to item (ii), we find by direct computation that $\|X\|_B = \|\m\|_1 \|\n\|_1.$  We define 2 positive unit  vectors $\xi$  and $\eta$ to be used in a bilinear factorization  of $X$ by  
 $$\eta_i := \frac{|\m_i|^{(1/2)}}{\|\m\|_1^{(1/2)} } , \quad  \xi_j := \frac{|\n_j|^{(1/2)}}{\|\n\|_1^{(1/2)} }.$$ 
The  matrix $ B$ in the factorization $X = \Delta_m(\eta)B\Delta_n(\xi)$ is then given by $$B_{(i,j)} =  \|\m\|_1^{(1/2)} \|\n \|_1^{(1/2)} 
\mathrm{sign}(\m_i)|\m_i|^{(1/2)} \mathrm{sign}(\n_j)|\n_j|^{(1/2)}.$$  Then, since the rank of $B$ is one,  $\|B\|_\infty = \|B\|_2 = \|\m\|_1\|\n\|_1$ and  we see from Theorem 1.3 item (ii) of \cite{C4}  that $$\|X\|_{cbB} \leq \|B\|_\infty = \|\m\|_1\|\n\|_1 =  \|X\|_B \leq \|X\|_{cbB},$$ so the norms are cross norms as claimed. 
By the definition of a minimal cross norm, it follows that $\|.\|_B$ is the minimal cross norm  $\vee(1,1) $ so the result $\|.\|_{cbB} \leq K_G^\bc \|.\|_{\vee(1,1)} $ follows from Theorem 2.4 of \cite{C3}.  The duality statement of item (ii) follows from equation (3.3) of \cite{C3} and item (ii) follows.  

With respect to item (iii), the equations $X_{(i,j)} = \m_i\n_j$ describes $ X $ as a product of a column matrix  $\m_|$  by a row matrix  $\n_{-},$ so by Theorem 1.3 item (iii)  of \cite{C4} we have  $\|X\|_S \leq \|\m\|_\infty \|\n\|_\infty.$ On the other hand a certain matrix unit will give the opposite inequality, and item (iii) follows together with a factorization $X = \m_| \n_{-}$ for which $\|\m_|\|_\infty \|\n_{-}\|_\infty  = \|X\|_S.$ The duality statement follows from that of item(ii), and the inequality $\|.\|_{\wedge(\infty, \infty)} \leq K_G^\bc\|.\|_S$ is the dual of the one presented in item (ii), so item (iii) follows.      

In the case of item (iv) we remark that $\|.\|_T$ is conjugate dual to the cross norm $\|.\|_{cbF}, $ so it is a cross norm, and then by the duality $\|X\|_T = \|\m\|_2\|\n\|_\infty.$  The concrete factorization may be obtained as follows. We consider again $X$ as a product of a column matrix  $\m_| $ and a row matrix $\n_{-}$. 
We have  $\|\Delta_m(\m)\|_2 = \|\m\|_2$  and for $L:=(\Omega_m)_{-}$ we have $\|L\|_c = 1$ and for $R = \n_{-}$ we have $\|R\|_c = \|\n\|_\infty. $ Then $X = \Delta_m(\m) L^*R$ so by item (iv) in Theorem 1.3 of \cite{C4} we get  $\|X\|_{T} \leq \|\m\|_2\|\n\|_\infty = \|X\|_T. $ The statement $\|.\|_{\wedge(2, \infty)} \leq \|.\|_T$ follows also from the conjugate  duality between the  cbF-norm  and the T-norm.  

The Haagerup tensor product is given by a cross norm $\|.\|_h$ on the tensor product of the two abelian C*-algebras  $(\ca_m, \|.\|_\infty)$ and $(\ca_n, \|.\|_\infty) $ in the following way for a matrix $X$ in $M_{(m,n)}(\bc).$  
 \begin{align*}\|X\|_h: = \inf \{ &  \sqrt{ \|\sum_k  a_k^*a_k\|  \|\sum_k b_k^*b_k\|}\, : \\ & a_k \in \ca_m, \, b_k \in \ca_n, \,  X = \sum_k \f(\overline{a_k} \otimes b_k)\,\}.
 \end{align*} 
If you are given an expression $X = \f(\sum_{k=1}^l  \overline{a_k} \otimes b_k) = \sum_k \overline{(a_k)}_|(b_k)_{-}, $ then you may construct an $l \times n $ matrix $R $ by $R_{(k,j)} := b_k(j)$ and an $l \times m $ matrix $L$ by $L_{(k,i)} := a_k(i).$ Then we find  that $L^*R = X, $ $\|R\|_c = \| \sum_k a_k^*a_k\|^{(1/2)}$ and $\|L\|_c = \|\sum_k b_k^*b_k\|^{(1/2)} $ and based on Theorem 1.3 item (iii) of \cite{C3} we have obtained a proof of the identity  $\|.\|_S = \|.\|_h,$ which is a part of item (v). The equation (3.9) of \cite{Pi3} and Theorem 1.3 item (iii) of \cite{C4} imply that the norm $\|.\|_H $ of \cite{Pi3} equals $\|.\|_S$ and item (v) follows. 

The equality  $\|.\|_H= \|.\|_S$  of item (v) implies by duality that $\|.\|_{H'} = \|.\|_{cbB}.$ The norm $\g_2$ from equation (2.2) of \cite{Pi3} is known to be equal to $\|.\|_{H'},$  and the entire theorem follows. 
\end{proof} 

Following Pisier's survey \cite{Pi3} we see that  all of the norms except  the norm $\|.\|_T$ are explicitly present in that survey, and as dual to one of these the $T-$norm is implicitly mentioned too.  The duality results mentioned in the items (ii) and (iii) was well known and used by Grothendieck. The real  new thing here  is - in our opinion -  that some of the norms now have a characterization as completely bounded norms of some linear or bilinear operators between operator spaces, and the theory of completely bounded linear and multi-linear maps yields concrete optimal factorizations of the operators under investigation.

\section{ The Haagerup factorization}  \label{Ha}

Haagerup gives in an appendix to \cite{Ha1} a proof of the little Grothendieck inequality for non commutative C*-algebras. That proof may also be applied to the case of the  finite dimensional commutative C*-algebra $\ca_n$ and then to a linear  operator $F_X$ from $(\bc^n, \|.\|_\infty)$ to $(\bc^m, \|.\|_2).$ In this way Haagerup provides  an elementary way to  obtain a concrete  factorization of a complex $m \times n$ matrix $X$  which qualitatively  is close to the optimal cbF-factorization of $X.$ See Theorem 3.1 item(i) of \cite{C3}.  Haagerup's construction  gives the best upper bound in the non commutative setting, but in the abelian case, which we study here, his  method gives the   upper bound 2 for for the little Grothendieck constant $k_G^\bc,$ and this aspect is not impressing, since we know that $k_G^\bc = 4/\pi < 1.274.$
The impressing thing is that his proof is constructive and only uses elementary analysis. Furthermore the construction shows that if we only look at matrices with non-negative  entries, then, in that world, $k_G^\bc = 1$ and Haagerup's factorization is the optimal one here. There is a possibility that Haagerup's factorization is the optimal one even if we get a bad upper bound for $k_G^\bc,$ the reason being that $\ca_n$ is abelian and in this case  the inequality \eqref{sa} below might have an extension to  a version of the inequality \eqref{nsa} with the constant $k_G^\bc $ instead of 2 ?  Even if Haagerup's construction does not give the optimal cbF-factorization in general,  it still  raises some new questions in this setting,  and it gives an elementary proof of Grothendieck's little inequality, although  it does not provide the right constant. 

 It is quite easy to describe Haagerup's factorization, and we will do that by following the first arguments in his proof of the non commutative little Grothendieck inequality. To this end, let $\Omega $ denote the vector in $\bc^n$  consisting of 1's only, and let  $X$ be a complex $m \times n $ matrix such that $\|X\|_F = 1.$ Then there exists a unitary $u$ in $\ca_n $ such that 
\begin{equation} \label{u} 1 = \|X\Delta_n(u)\Omega\|_2 , \end{equation} and we can define  a functional  $\phi$ on $\ca_n$ of norm  $1$ by 
$$ \phi(a) := \langle X \Delta_n(ua) \Omega, X\Delta_n(u) \Omega \rangle = \langle \Delta_n(a) \Omega , \Delta_n(u)^*X^*X\Delta_n(u) \Omega\rangle.$$
This functional of norm 1 on $\ca_n = C(\{1, \dots, n\}, \bc)$  takes the value 1 at the unit $I_n,$ so it is a state - i. e. given as the integral with respect to a probability measure on the set $\{1, \dots, n\}$. 
Then the vector $\l$ in $\bc^n$ defined by  $$\l = \Delta_n(u)^*X^*X\Delta_n(u)\Omega $$ has non negative entries with sum 1, and these real numbers are  the masses of the points with respect to the mentioned  measure. We can then define a unit  vector $\xi$  with non-negative entries   in $(\bc^n, \|.\|_2)$  by the definitions $\xi_j = \sqrt{\l_j}.$ 
If we follow Haagerups proof, we find  that   
\begin{align} 
\forall a =a^* \in \ca_n :&   \label{sa} \\ \notag \quad  \|X\Delta_n(u) \Delta_n(a) \Omega\|^2_2 &= \|F_X(ua)\|^2_2  \leq  \|\Delta_n(a)\xi\|^2_2 = \|\Delta_n(\xi) \Delta_n(a) \Omega\|^2_2\\ 
\forall a  \in \ca_n  :& \label{nsa} \\ \notag \quad \|X\Delta_n(u) \Delta_n(a) \Omega\|^2_2  &  \leq 2 \|\Delta_n(a)\xi\|^2_2 = 2\|\Delta_n(\xi) \Delta_n(a) \Omega\|^2_2. 
\end{align}

Equation \eqref{nsa} shows that there exists a complex  $m \times n $ matrix $Z$ such that $\|Z\|_\infty \leq \sqrt{2},$  and $X\Delta_n(u) = Z\Delta_n(\xi) $ and $X = (Z\Delta_n(u^*)) \Delta_n(\xi).$ This gives the rather bad upper bound $k_G^\bc \leq 2$ for the general  cbF-factorization, but a quite explicit construction. We have made some simple experiments with this vector in the mathematical tools package named {\em Maple,} but our knowledge of  the powerful tools  of Maple is limited, so no new information showed up, except for the case of matrices with non-negative entries. This case is much easier to deal with, because here the optimal unitary $u$ which is a part of Haagerup's construction is simply the unit $I_n$ of $\ca_n.$ The experiments for matrices with non-negative entries  indicated that for such matrices we have$$ \|X\Delta_n(\xi)^{-1}\|_\infty = \|X\|_F. $$  If this experimental result is true then it may be combined  with the results from \cite{C3} to get  $$\|X\|_{cbF} \leq \|X\Delta_n(\xi)^{-1}\|_\infty = \|X\|_F \leq \|X\|_{cbF},$$ which, if true,  shows  that the Haagerup vector is the optimal one for the cbF-factorization in this case.  The following easy proposition and our reformulation of Grothendieck's inequality as a linear program in Theorem 3.5 of  \cite{C4} tell that the results of the  experiments are based on a mathematical theorem. This is yet another example which demonstrates  that the points of view presented in the  book \cite{EJ} by Eilers and Johansen  
are very fruitful.

\begin{proposition} \label{NxF} 
Let $X$ be a complex $m \times n$ matrix, then the 1-norms  $\|_iX\|_1 $ of the rows in $X$ satisfy the inequality  \begin{equation} \label{PosF}  \|X\|_F \leq \big(\sum_i \|_iX\|_1^2\big)^{(1/2)}.\end{equation}  If the entries in $X$ are all non negative then \eqref{PosF} becomes an equality.

\begin{proof}
Let $z$ in $\bc^n$ be given with $\|z\|_\infty \leq 1,$ then \begin{align*}
\|F_X(z)\|_2^2 & = \sum_i |\sum_j X_{(i,j)}z_j|^2 \\
& \leq \sum_i ( \sum_j |X_{(i,j)}| )^2 \\
&= \sum_i \|_iX\|_1^2,
\end{align*} and the proposition follows.
\end{proof}
\end{proposition}  
 
\begin{theorem}\label{PosHa}
Let $X$ be a real $m \times n$ matrix with non negative entries then $\|X\|_{cbF} = \|X\|_F = \big(\sum_i \|_iX\|_1^2\big)^{(1/2)}, $ and Haagerup's vector becomes the optimal positive unit vector from the cbF-factorization of $X.$ 
\end{theorem} 

\begin{proof}
We define  a positive matrix $P = X^*X$ and a vector $\l$ in $\br^n$ with non negative entries via the formula $$ \l_j := \sum_s \sum_t x_{(s,j)}x_{(s,t)},$$  and we will show that $P \leq \Delta_n(\l).$  Then let $\g$ be a vector in $\bc^n$ and we get \begin{align} \label{d}
\langle P\g, \g\rangle & = \|X\g\|^2 \\ \notag & = \sum_s| \sum_j X_{(s,j)}\g_j|^2  \text{ by Cauchy-Schwarz's Inequality  } \\ \notag  & \leq \sum_s\bigg(\sum_j X_{(s,j)} |\g_j|^2 \bigg)\bigg(\sum_t X_{(s,t)}\bigg) \\ \notag & = 
\sum_j \l_j |\g_j|^2\\ \notag & = \langle \Delta_n(\l) \g, \g\rangle.
\end{align}
We see by the proposition, the definition of $\l$  that $\sum_j \l_j = \|X\|_F^2, $  then by Theorem 3.5 of \cite{C4} and the computations above we get  that $\|P\|_{cbB} \leq \sum_j \l_j.$ By  Theorem 3.2 of \cite{C4}  we know that $\|P\|_{cbB} = \|X\|_{cbF}^2,$ so $$\|X\|_{cbF}^2= \|P\|_{cbB} \leq \sum_j \l_j = \|X\|_F^2\leq \|X\|_{cbF}^2.$$ By Theorem 3.5 of \cite{C4} we know that the optimal positive unit  vector $\xi$ for the cbF-factorization is given by the formula $\xi_j = \sqrt{\l_j/ (\|X\|_F^2)} $ and that is Haagerup's vector, so the corollary follows.  
\end{proof}
 The results  above also give an upper bound -  which is easy to compute - for the cbB-norm of a positive matrix.  

 \begin{corollary}
 Let $P$ be a positive $n \times n $ complex matrix and $X$ a complex  $m \times n$ matrix with $X^*X = P$ then \begin{align*} \mathrm{Tr}_n(P) \leq \|P\|_{cbB} &\leq \sum_s \sum_j \sum_t |X_{(s,j)}||X_{(s,t)}| 
 \leq \big( \mathrm{Tr}_n\big( \mathrm{diag}(P)^{(1/2}\big)\big)^2,
 \end{align*} with equality in the second inequality if all the entries $X_{(i,j)}$ are non negative. 
 \end{corollary} 
\begin{proof}
We return to \eqref{d}, and we find that for the non negative   $\l    $ in  $\br^n$  defined by 
$$ \l_j := \sum_s \sum_t |X_{(s,j)}||X_{(s,t)}|$$
we will get  $P \leq \Delta_n(\l)$ and then by Theorem 3.5 of \cite{C4}, $\|P\|_{cbB} \leq \mathrm{Tr}_n(\Delta_n(\l)), $ which is the stated sum. To get the next inequality we  remind you that the column norm of the $j'$th column $\|X_j\|$  equals $P_{(j,j)}^{(1/2)}.$ We can then continue with the inequality we already have established and use the  Cauchy-Schwarz Inequality once more  to get 
\begin{align*}
\|P\|_{cbB} & \leq \sum_s \sum_j \sum_t |X_{(s,j)}||X_{(s,t)}| \\ & \leq \sum_j \big( \sum_u |X_{(u,j)}|^2 \big)^{(1/2)} \sum_t \big( \sum_v |X_{(v,t)}|^2 \big)^{(1/2)} \\& =  (\sum_j\|X_j\|)^2 \\ &= \big( \mathrm{Tr}_n( \mathrm{diag}(P)^{(1/2})\big)^2,
\end{align*} and the corollary follows.
\end{proof}

It seems natural to ask if the results above can be extended to matrices with real entries only? Since it is known that  $k_G^\br > 1,$ this is not possible, but why not.  One  reason, we anticipate, is that for some matrices with real entries  the maximal value of $\|F_X(u)\|_2$ must be obtained in a unitary $u = (u_1, \dots, u_n)$ where the entries are not all real. 
In the general case, it might still be that Haagerup's vector continues to be the one from the optimal cbF-factorization, but we think that this is  not the case. We do not have a definite result which proves this, but we have made some observations, which indicate this and also give some more information on the problems in finding the optimal vector $\xi$ for the cbF-factorization.  

Suppose we are given a complex $m \times n$ matrix $X$ with optimal unitary $u$ in $\ca_n$ and a Haagerup vector $\xi$ such that all $\xi_j >0, $  then we have  the following equation $$\frac{1}{\|X\|_F^2}  \bigg(\Delta_n(\xi)^{-1} \Delta_n(u^*) X^*X \Delta_n(u) \Delta_n(\xi)^{-1}\bigg)  \xi = \Delta_n(\xi)^{-1} \l = \xi, $$
 so for $\g := \Delta_n(u) \xi $ get
 $$ \frac{1}{\|X\|_F^2} 
 \bigg(\Delta_n(\xi)^{-1} X^*X  \Delta_n(\xi)^{-1}\bigg)  \g = \g$$ and $\g$ is an   eigenvector corresponding to the eigenvalue $\|X\|_F^2$  for the matrix 
  $\Delta_n(\xi)^{-1}  X^*X  \Delta_n(\xi)^{-1} .  $  By the factorization result in Theorem 3.1 item (i) of \cite{C3} we know that $\|\Delta_n(\xi)^{-1}  X^*X  \Delta_n(\xi)^{-1}\|_\infty \geq \|X\|_{cbF}^2$ 
  and with equality exactly if $\xi$ is the positive unit vector from the cbF-factorization of $X.$ 
 For a general matrix $X$ we can then use Haagerup's recipe to construct a positive unit vector $\xi$ such that  $\|X\|_F^2 $ 
 is an eigenvalue for the matrix $\Delta_n(\xi)^{-1}  X^*X  \Delta_n(\xi)^{-1}$  and
  $$\|X\|_F^2 \leq \|X\|_{cbF}^2 \leq \|\Delta_n(\xi)^{-1}  X^*X  \Delta_n(\xi)^{-1}\|_\infty  \leq 2 \|X\|_F^2.$$
  If the recipe actually gave the optimal cbF-factorization from Theorem 3.1 of \cite{C3} then  both $\|X\|_F^2$ and $\|X\|_{cbF}^2$ will  be eigenvalues for  $\Delta_n(\xi)^{-1}  X^*X  \Delta_n(\xi)^{-1},$  and that may be too much to ask for ?
  We have obtained a way to reformulate this problem, and in order to describe this we return to the the question on computing the cbB-norm of $X^*X.$  Then  let $X^*X = \Delta_n(\eta)B\Delta_n(\eta) $ be its cbB-factorization and we can  define a probability distribution $\mu$ on $\{1, \dots, n\}$ by $\mu_j = \eta_j^2.$ Suppose all $\eta_j > 0,$ 
  then from Theorem 2.4 of \cite{C3} we know that $\|X^*X\|_{cbB} $ is an eigenvalue of $B$ and then the following determinants both vanish, $$ \mathrm{det}\big(X^*X - \|X^*X\|_{cbB}\Delta_n(\mu) \big) =0 \text{ and } \mathrm{det}(\big(X^*X - \|X\|_F^2 \Delta_n(\l) \big) =0,$$  and we have obtained the following proposition. 
 
 \begin{proposition} \label{prop} Let $X$ be a complex $m \times n $ matrix, $\l$ the probability distribution associated to Haagerup's vector and $\mu$ the probability distribution associated to the cbB-factorization of $X^*X.$ If all $\l_j >0$ and all $\m_j>0$ then both $\|X\|^2_F\Delta_n(\l)$ and $\|X^*X\|_{cbB}\Delta_n(\m)$ belong to the set of positive diagonal matrices  $D$ 
 which satisfy $\mathrm{det}\big(X^*X - D \big) = 0.$ 
 \end{proposition}  
We do not know if the equation $\mathrm{det}\big(X^*X - D \big) = 0$ of Proposition \ref{prop} has been the  studied in the literature, but the structure of its set of solutions might offer some new insights into the questions we are looking at here. 

\section{From the little inequality to the Grothendieck inequality. The analytic approach.} \label{constants} 

In Theorem 3.2 of \cite{C4}
 we showed that for any complex $m \times n$ matrix $X$ we have $\|X\|_F^2 = \|X^*X\|_{B} $ and $\|X\|_{cbF}^2 = \|X^*X\|_{cbB}.$ 
This shows that $$\frac{\|X^*X\|_{cbB}}{\|X^*X\|_B} =    \bigg(\frac{\|X\|_{cbF}}{\|X\|_F}\bigg)^2,$$  hence  the following equation is valid, and  $k_G^\bc$ is the smallest possible constant for this inequality.  \begin{equation} \label{K=kPos}  \forall P  \in M_n(\bc)_+ : \quad  \|P\|_{cbB} \leq k_G^\bc \|P\|_B.
 \end{equation}
Then  in the world of positive square matrices the Grothendieck inequality has the constant $k_G^\bc.$ In particular this implies the well known inequality $k_G^\bc \leq K_G^\bc.$ To try to understand how $K_G^\bc$ depends on $k_G^\bc$ for non positive matrices, we will in Theorem \ref{Fact}  present a factorization result,  which improves some of the non self-adjoint aspects of Theorem 3.2 of \cite{C4}.  Unfortunately we were not able to make direct use  of   this result in the way we hoped for, but some steps in its proof are actually used in the proof of Theorem \ref{KG}, where we show that the general Grothendieck's inequality does follow from the inequality for positive matrices with the constant  $K_G $ satisfying  $ k_G^\bc < K_G^\bc \leq k_G^\bc/(2- k_G^\bc) < 1.752.$   

\begin{theorem} \label{Fact} 
Let $X$ be a complex $m \times n $ matrix and let the cbB-factorization of  $X$ be denoted $X  = \Delta_m(\eta) B \Delta_n(\xi).$ Let $B= WP $ be the polar decomposition of $B, $ then the complex  $n \times n$ matrix  $C$ and the complex $ m \times n$ matrix  $D$ given as $ C : = P^{(1/2)}  \Delta_n(\xi )$ and $D:= \Delta_m(\eta) WP^{(1/2)}$ satisfy \begin{itemize} \item[(i)] $X = DC.$ 
\item[(ii)] $\|D^*\|_{cbF} = \|X\|_{cbB}^{(1/2)} $ and $D^*= (P^{(1/2)}W^*)\Delta_m(\eta)$ is the  cbF-factorization of $D^*.$ 
\item[(iii)] $\|C\|_{cbF} = \|X\|_{cbB}^{(1/2)}$ and $C= P^{(1/2)}\Delta_n(\xi)$ is the  cbF-factorization of $C.$ 
\end{itemize} 
\end{theorem}
\begin{proof}
Suppose for simplicity that $\|X\|_{cbB} =1.$ By the theorems 1.8 and  1.3 item (iii) in \cite{C4}, choose a matrix $L $ in $M_{(n,m )}(\bc)$ with  $\|L\|_c = 1$ and a matrix $R$  in $M_n(\bc)$ with $\|R\|_c = 1 $ such that Tr$_n (R^*L X ) = 1.$ Then \begin{align*}
1 &= \mathrm{Tr}_n\big((L\Delta_m(\eta ) W P^{(1/2)})( P^{(1/2)} \Delta_n(\xi)R^*)\big) \\ 
&= \langle  L\Delta_m(\eta ) W P^{(1/2)}, R\Delta_n(\xi )  P^{(1/2)}  \rangle.
\end{align*}
We can see - as in equation (2.1) from \cite{C4} - that since $\|B\|_\infty =1 $ we have  $\|L\Delta_m(\eta ) W P^{(1/2)}\|_2 \leq 1,$ and $\|R\Delta_n(\xi )  P^{(1/2)} \|_2 \leq 1 $ then  both inequalities are equalities and we can define an $ n \times n $ complex matrix $T$ with $\|T\|_2 = 1$ by the equations  
\begin{equation}
T:= P^{(1/2)} W^* \Delta_m(\eta) L^*  = P^{(1/2)} \Delta_n(\xi)R^* .
\end{equation}
Since $\|P\|_\infty = 1 $ we get from Theorem 1.3 item (i) of \cite{C4} that $$\|\Delta_m(\eta)W P W^* \Delta_m(\eta)\|_{cbB} \leq 1 \text{ and } \|\Delta_n(\xi)P\Delta_n(\xi) \|_{cbB} \leq 1.$$ Similarly $\|L\|_c \leq 1 $ and $\|R\|_c \leq 1 $ imply that $\|L^*L\|_S \leq 1 $ and \newline $\|R^*R \|_S \le 1. $ Then \begin{align}
\notag 1 &= \mathrm{Tr}_n(T^*T) \\&= \mathrm{Tr}_m\big( (\Delta_m(\eta)W PW^* \Delta_m(\eta))(L^*L)\big) \label{L*L}\\& \notag \leq \|\Delta_m(\eta)W PW^* \Delta_m(\eta)\|_{cbB} \|L^*L\|_S \leq 1 \\ & \notag \text{ and also } \\ \notag  1  &= \mathrm{Tr}_n(T^*T) \\&= \mathrm{Tr}_n\big( (\Delta_n(\xi)P\Delta_n(\xi))(R^*R )\big) \label{R*R} \\& \notag\leq \|\Delta_n(\xi) P\Delta_n(\xi)\|_{cbB} \|R^*R\|_S \leq 1 .
\end{align}
The theorem then follows from Theorem 3.6 of \cite{C3}. 
\end{proof}
We will now apply the theorem above and  parts of its    proof to get an upper bound for $K_G^\bc$ expressed in terms of $k_G^\bc.$ 

\begin{theorem} \label{KG}
Grothendieck's complex constants satisfy $$ k_G^\bc \leq K_G^\bc \leq \frac{k_G^\bc}{2 - k_G^\bc}.$$
\end{theorem} 
\begin{proof}
As mentioned in front of Theorem \ref{Fact}, that theorem, when applied to poitive matrices, implies that $k_G^\bc \leq K_G^\bc.$ Let   $X$ and $Y$  be  complex $m \times n$ matrices such that $\|X\|_{cbB} = 1,$ $\|Y\|_S = 1$ and Tr$_n(Y^*X) = 1. $ We will also assume that no row and no column in $X$ vanishes. Let the cbB decomposition of $X$ be denoted as  $X = \Delta_m(\eta)B\Delta_n(\xi )$ and an elementary Schur decomposition of $Y $ be given as $Y =  L^*R$ with $L$ in $M_m(\bc),$  $ R$ in $M_{(m,n)}(\bc),$ and  $\|L\|_c = \|R\|_c=1.$  The non-vanishing of rows and columns in $X$ imply that all $\eta_i >0$ and all $\xi_j > 0.$
We will now construct a positive matrix $P$ in $M_{(m+n)}(\bc),$ to which we will apply equation\eqref{K=kPos},  and we will also define  a positive matrix $Q$ in $M_{(m+n)}(\bc)$ with Schur multiplier norm 1 and other nice properties. Let $\g$ be the positive unit vector in $\bc^{(m+n)} $ given as $\g = 2^{-(1/2)}( \eta_1, \dots, \eta_m, \xi_1, \dots , \xi_n),$ then   we define $P$ and $Q$ by

\begin{align} \notag P :&= \begin{pmatrix}
\Delta_m(\eta)^2 & X\\  X^* & \Delta_n(\xi)^2
\end{pmatrix}\\ \notag  &= \begin{pmatrix} \Delta_m(\eta) &0\\0 & \Delta_n(\xi )
\end{pmatrix} \begin{pmatrix}I_m & B\\ B^* & I_n\end{pmatrix} \begin{pmatrix}\Delta_m(\eta) & 0\\ 0&\Delta_n(\xi)
\end{pmatrix}\\ \notag
&= 2\Delta_{(m+n)}(\g)
 \begin{pmatrix}I_m & B\\ B^* & I_n\end{pmatrix} \Delta_{(m+n)}(\g)\\ \label{Q}
 Q :&= \begin{pmatrix}
L^*L  & L^*R\\ R^*L & R^*R\end{pmatrix}\\ \notag   &= \begin{pmatrix} L^*\\ R^*
\end{pmatrix}\begin{pmatrix} L & R
\end{pmatrix}.
\end{align}
The operator $B$ has operator norm 1 since 
$\|X\|_{cbB} =1. $ Then the matrix $\begin{pmatrix}
I_m & B \\B^* & I_n\end{pmatrix}$  is positive and of operator norm 2. This implies that $P$ is positive and $\|P\|_{cbB} \leq 4.$ The operator $Q$ is clearly positive and diag$(Q) \leq I_{(m+n)}$ since $\|L\|_c =1 $ and $\|R\|_c =1.$ Since $\|L\|_c = 1$ we know  by Schur's result \cite{Sc}  that  $\|Q\|_S = 1.$ The last decomposition of $Q$ 
also shows  that $\|Q\|_S =  1 $ since the column norm of the matrix $(L\, \, \, R) $ is 1. As in the proof of Theorem \ref{Fact}, the  equality Tr$_n(Y^*X) = 1,$ implies that the analogies to the equations \eqref{L*L} and \eqref{R*R}  hold here, too. Since $\|B\|_\infty =1$ we get from \eqref{L*L} that Tr$_m(\Delta_m(\eta)L^*L\Delta_m(\eta) = 1.$ 
Similarly  \eqref{R*R} implies  that Tr$_n(\Delta_n(\xi)R^*R\Delta_n(\xi) = 1,$ and combined we get that  Tr$_{(m+n)} ( QP) =  4,$ so $\|P\|_{cbB} = 4. $ We will then  find the value of $\|P\|_B.$ It is not hard to see that we get $\|P\|_B = 2 + 2\|X\|_B.$ We already mentioned that Theorem \ref{Fact} of this article or more explicit Theorem 3.2 of \cite{C4} implies that for the positive matrix $P$ we have $\|P\|_{cbB} \leq k_G^\bc\|P\|_B $ and then \begin{align} \notag
4 &\leq k_G^\bc (2 + 2\|X\|_B) \\ \notag
\|X\|_B & \geq \frac{2-k_G^\bc}{k_G^\bc} .
\end{align} 
This holds for any $X$ with $\|X\|_{cbB} =1, $ so for a general $X$ we get 
\begin{equation} \label{kK} 
\|X\|_{cbB}  \leq \frac{k_G^\bc}{2-k_G^\bc} \|X\|_B,\end{equation} and the theorem follows. 
\end{proof}

\section{A geometrical characterization of $k_G^\bc$  and a geometrical proof of Grothendieck's inequality} 
The result \eqref{K=kPos} for positive   $P,$  that  $\|P\|_{cbB} \leq k_G^\bc \|P\|_B$ and the fact that $k_G^\bc$ is the smallest possible constant, which satisfies this inequality, must clearly be based on geometrical properties  of some unit balls with respect to some norms, but it is not obvious what the relevant relations  may be. We are missing a lot in our understanding of the geometrical aspects, but we have found a relation between some convex sets which reflects a geometrical property of the positive part of the unit ball of the Schur multipliers. Our  result gives a version of Grothendieck's   geometrical formulation of his inequality as presented  in Theorem 1.1. We will obtain the constant $k_G^\bc/(2- k_G^\bc)$ in stead of $K_G^\bc.$ This approach also leads to a geometrical characterization of the value of $k_G^\bc.$ In this situation we look at positive   Schur multipliers and try to write them as linear combinations of positive  rank one multipliers in a fashion which is close to the optimal one from Theorem 1.1. If we look at a positive $n \times n$ matrix $P$ with $\|P\|_S =1,$ the natural desire would then be to write it 
as a positive linear combination of rank one positive  matrices in  the form  $$P = \sum_k \a_k (\bar \xi_k)_| (\xi_k)_{-}, \, \, \a_k \geq 0,\, \sum_k \a_k \leq K_G^\bc,\, \xi_k  \in \bc^n,\,\|\xi_k\|_\infty \leq 1.$$  This hope is a fake dream, as far as we can see, and the extension of Theorem 1.1 to a neater result for positive matrices is much more complicated. The complications can be done with, but at a cost of a new point of view at the geometrical problem. The following Theorem \ref{Chark}  shows that  $k_G^\bc$ may be seen as a constant, which expresses a relation between certain compact convex sets. This point of view reproduces \eqref{kK} as an easy corollary.  
  
  \begin{definition}
  \begin{itemize}
\item[(i)]    The compact convex subset $\cq_n$ of $M_n(\bc)$ is defined by $$ \cq_n  = \{ X \in M_n(\bc) \,: \, X \geq 0 \text{ and } \mathrm{diag}(X) = I_n\}.$$
 \item[(ii)] The compact convex subset $\car_n$ of $M_n(\bc)$ is defined as the closed convex hull of the set $\{ ( u^*)_| (u)_{-}\, : \, u \text{ unitary in } \ca_n\}.$ 
 \end{itemize}
  \end{definition}
It is worth to keep in mind that $\car_n$ is the closed convex hull of the rank one matrices in $\cq_n$ so it is clearly a subset of $\cq_n.$ 
The geometrical statement in this section is the following Theorem \ref{geo}, but we will need a little lemma first.   
 \begin{lemma}  \label{QRnorms}  

Let $P$ be a positive $n \times n$ matrix. Then there exists an $R$ in $\car_n$ and a $Q$ in $\cq_n$ such that $\mathrm{Tr}_n(PR) = \|P\|_B$ and $\mathrm{Tr}_n(PQ) = \|P\|_{cbB}.$  \end{lemma} 
\begin{proof}
By Theorem 3.2 of \cite{C4} we know $\|P\|_B^{(1/2)} = \|P^{(1/2)}\|_F.$ The latter norm is the operator norm for $P^{(1/2)}$ as an operator from the C*-algebra $\ca_n$ into the Hilbert space $\bc^n.$ The extreme points in the unit ball of $\ca_n$ are the unitaries 
so there exists a unitary $u$ in $\ca_n$  and  an $R = u_|\bar u_{-}$ in $\car_n$ such that $$\mathrm{Tr}_n(PR) = \|P^{(1/2)}u_| \|^2 = \|P^{(1/2)}\|_F^2 = \|P\|_B.$$   
By Theorem 1.8 of \cite{C4} we know that there exists   a self-adjoint $Q_0$ in $M_n(\bc)$ such that $\|Q_0\|_S = 1 $ and $\mathrm{Tr}_n(PQ_0) =\|P\|_{cbB}.$ The Proposition 2.4 of \cite{C4} shows that $Q_0 $ has a factorization $Q_0 = GSG$  such that $G $ is positive with diag$(G^2) \leq I_n$ and $ S$ is a self-adjoint partial isometry, which in this case means an orthogonal  projection. 
We may then define  a positive diagonal operator $D $ by $D:= I_n - \mathrm{diag}(G^2)$ and an element $Q$ in $\cq_n$ by $Q := G^2 +D,$ so the diagonal of $Q$ is $I_n$ and $\|Q\|_S = 1.$ Hence $$\|P\|_{cbB} = \mathrm{Tr}_n(PQ_0) \leq \mathrm{Tr}_n(PG^2) \leq  \mathrm{Tr}_n(PQ) \leq \|P\|_{cbB}\|Q\|_S = \|P\|_{cbB}$$ and the lemma follows. 
\end{proof}

 \begin{proposition} \label{geo}  
 Let $Q$ be a matrix in $\cq_n$ then there exist matrices   $R$ in $\car_n$ and $ P$ in $M_n(\bc)_+$ such that $  Q =  k_G^\bc R - P.$     
 \end{proposition}

\begin{proof}
We will use the duality between the norms $\|.\|_{cbB} $ and $\|.\|_S.$ It follows from the computations in \cite{C4} that this duality also holds if we restrict to the real subspace $M_n(\bc)_{sa} $ consisting of the self-adjoint complex $n \times n$  matrices. The convex cone  $M_n(\bc)_+$ consists of the positive matrices in $M_n(\bc).$ We will work inside the real vector space consisting of the  self-adjoint matrices and here we define the polar of a set $\cas $ contained in $M_n(\bc)_{sa} $ by $$ \cas^\circ := \{X \in M_n(\bc)_{sa} \, : \, \forall  S \in \cas : \mathrm{Tr}_n(SX) \leq 1\}.$$ 
We define two convex subsets $\cas $ and $\ct$ of $M_n(\bc)_{sa} $ by 
\begin{align}
\cas :&= \cq_n - M_n(\bc)_+ \\
\ct :&= \car_n - M_n(\bc)_+,
\end{align}
then we can remark that bot sets contain the zero matrix, and they are both closed since both $\car_n$ and$\cq_n$ are compact. Hence both sets equal their bipolars.  In order to compute their polars we need a little observation based on the theorems 1.3 plus 1.8 of \cite{C4}. 

Based on standard techniques and the lemma above we get  
\begin{align*}
\cas^\circ &= \{ X \in M_n(\bc)_+ \, : \, \|X\|_{cbB} \leq 1\}. \\
\ct^\circ &= \{ X \in M_n(\bc)_+ \, : \, \|X\|_{B} \leq 1\}.
\end{align*}
Equation \eqref{K=kPos}  shows that $\ct^\circ \subseteq k_G^\bc \cas^\circ,$ and then  by the bipolar theorem 
\begin{equation} \label{Subset} 
\cq_n \subseteq k_G^\bc(\car_n - M_n(\bc)_+) = k_G^\bc \car_n  - M_n(\bc)_+,  
\end{equation} and the proposition  follows.
\end{proof}
We will use  equation \eqref{Subset} to obtain 2 results. First the equation \eqref{Subset} characterizes $k_G^\bc$ in a geometrical way,  and secondly it can give a geometrical proof of Theorem 1.1, although with a constant larger than $K_G^\bc.$ 

\begin{theorem}  \label{Chark}
The constant $k_G^\bc$ is the smallest positive real  $\b$ such that for any natural number $n$ and any real $\a$ with $\a \ge \b$  we have 
\begin{equation} \label{eqn} 
\cq_n \subseteq \a \car_n   - M_n(\bc)_+.
\end{equation}
\end{theorem}

\begin{proof}
We know from \eqref{Subset} that \eqref{eqn} holds for $\a = k_G^\bc,$ so we will first remark that if \eqref{eqn} holds for an $\a > 0 $ then it holds for all $\g \geq \a.$ 
To see this we just make the following rearrangement for a $Q$ in $\cq_n.$ By the assumption there exists an $R$ in $\car_n$ and a $P$ in $M_n(\bc)_+$ such that $Q = \a R - P = \g R - \big(P+ (\g-\a)R\big), $ and the claim follows.
 We can then define $\a_n $ as the infimum over all the possible $\a$'s which, for a given $n,$  satisfy \eqref{eqn} and because $\car_n$ is compact this $\a_n$  will also satisfy that equation. 
 Finally we define $\b$ as the supremum over all the $\a_n$'s.  By the equation  \eqref{Subset} we get $\b \leq k_G^\bc.$ On the other hand let $S$ be a positive $n \times n $ matrix then by Lemma \ref{QRnorms} there exists a $Q$ in $\cq_n$ such that Tr$_n(QS) = \|S\|_{cbB}. $ To this $Q$ we can find $R$ in $\car_n$ and $P$ in $M_n(\bc)_+$ such that $Q = \b R -P,$ and then \begin{equation}
 \|S\|_{cbB} = \mathrm{Tr}_n(QS) \leq \b \mathrm{Tr}_n(RS) \leq \b \|S\|_B. 
 \end{equation}
By the statements in front of the inequality \ref{K=kPos}, or at the beginning of this section we get $\b \geq k_G^\bc,$ and the theorem follows. 

The following proposition is an application of Theorem \ref{geo} which will serve to get a geometrical proof of Grothendieck's inequality for Schur multipliers.

\begin{proposition} \label{geo2}Suppose a natural number $n$ is given and let $Q$ be a matrix in $\cq_n$ then there exist matrices   $R_+$ and $ R_-$ in $\car_n $ such that 
$$  Q =  \frac{1}{2-k_G^\bc} R_+ - \frac{k_G^\bc -1}{2 - k_G^\bc} R_-.$$     
 \end{proposition}   

 To a matrix  $Q$ in $\cq_n$ there exist  by \eqref{Subset}   a matrix $R_1$ in $\car_n$ and a positive matrix $P_1$ such that \begin{equation} 
Q = k_G^\bc R_1 - P_1.
\end{equation}
From this we see that diag$(P_1) = (k_G^\bc -1)I_n,$ so $P_1 $ is in the set $(k_G^\bc -1) \cq_n,$ and we can use the equation \eqref{Subset} once more to obtain the existence of a matrix $R_2$ in $\car_n$ and a positive matrix $P_2$ such that 
\begin{align*}
P_1 &= k_G^\bc(k_G^\bc -1) R_2 - P_2, \text{ then } \\
Q & = k_G^\bc R_1 - k_G^\bc(k_G^\bc -1)R_2 + P_2\\
\mathrm{diag}(P_2)& = (k_G^\bc -1)^2 I_n.
\end{align*}
We can then continue by induction and since $ 0 < (k_G^\bc - 1) < 1 $ we can obtain convergent sums to describe $Q$ as $$ Q = \big( \sum_{k=1}^\infty k_G^\bc(k_G^\bc -1)^{2(k-1)} R_{(2k-1)}\big) - \big(\sum_{k=1}^\infty k_G^\bc(k_G^\bc -1)^{(2k-1)} R_{(2k)}\big)$$ and then by convexity and closedness of $\car_n,$ there exist matrices $R_+$ and $R_- $ in $\car_n$ such that \begin{equation} \label{Qdiff} 
 Q = \frac{1}{2-k_G^\bc}R_+ - \frac{k_G^\bc -1}{2-k_G^\bc}R_-,
\end{equation} and the proposition follows. 
\end{proof}

It is possible to obtain the Schur multiplier variant of Theorem \ref{KG} as a corollary to this proposition,  so we will present this proof too. First we define $\cv$ as the closed convex hull of all rank 1 matrices where all entries have numerical value equal to 1. This set may be described as the convex hull of  the matrices which are given as the product of a unitary $m-$column vector by a unitary $n-$row vector \begin{equation}
\label{T} \cv:= \overline{\mathrm{conv}}\big(\{ (u_|v_{-}) \,:\, u \text{ unitary in  } \ca_m, \, v \text{ unitary in } \ca_n\,\}\big).  
\end{equation} 
It is worth to remark that although $\cv$ and $\car_n$ seems to have similar definitions, they are anyway very different as sets, the reason beeing that in the definition of $\cv$ the variables $u$ and $v$ are independent. For instance  $0$ is in $\cv$  while for any $R$ in $\car_n$ we have diag$(R) = I_n.$ We can now state the corollary.  
\begin{theorem} \label{KGg}
Let $ X$ be a complex $m \times n$ matrix such that $\|X\|_S =1,$ then $X$ belongs to the set $k_G^\bc/(2-k_G^\bc) \cv.$
\end{theorem} 
\begin{proof}
By the proof of Proposition 2.6 of \cite{C4} we see that the condition $\|X\|_S = 1$ implies that there exists a natural number $k,$ a $k \times m $ matrix $L$ and a $k \times n$ matrix $R$ such that $X = L^*R$ and any column in both $L$ and $R$ is a unit vector. Then the matrix $Q$  in $M_{(m+n)}(\bc) ,$ which is defined in \eqref{Q} is positive and has diagonal equal to $I_{(m+n)},$ so we may apply the Proposition \ref{geo2}. If we look at the $(1,2)$ corner of the block matrix $Q$ we find that  $X = L^*R$ sits there and satisfies $$X \in \big(\frac{1}{2 - k_G^\bc} + \frac{k_G^\bc -1}{2-k_G^\bc} \big)\cv = \frac{k_G^\bc}{2 - \k_G^\bc}\cv,$$ so the theorem follows.
\end{proof}

The theorem raises the natural question if the equation 
\eqref{Qdiff} is valid with another  constant, different from $1/(2-k_G^\bc),$  such that the corollary would give Theorem 1.1 with the right constant $K_G^\bc.$ This is not possible, at least if one follows the obvious path to try. The following lines contain the analysis which leads to this conclusion.  

Remark that if for some  positive real $\a \geq 1 $ we have \begin{equation} \label{alpha}  \cq_n =\big( \a \car_n - (\a-1) \car_n\big) \cap M_n(\bc)_+ ,\end{equation}  then this relation holds with $\a$ replaced by any real $\b \geq \a.$ This follows because $\car_n $ is convex and it is seen in the following way. 
Let $Q$ in $\cq_n$ and $R_1, R_2$ be in $\car_n$ such that $Q = \a R_1 - (\a-1)R_2,$ then we define $R_3$ in $\car_n$ by $$ R_3:=  \frac{\b - \a }{\b -1}R_1 + \frac{\a-1}{\b-1}R_2,$$ and we get $Q = \b R_1 -(\b-1)R_3$. Since $\car_n$ is compact there exists a minimal $\a,$ say $\a_n,$ such that \eqref{alpha} is valid in $M_n(\bc)_+$ for $\a_n$ and all reals larger than $\a_n.$  
The question is then which relation is there between  $\a_s:= \sup\{\a_n\, : \, n \in \bn\}$ and $K_G^\bc $ or $k_G^\bc$ ? We have some partial answers.  By Theorem \ref{geo} we get $\a_s \leq 1/(2-k_G^\bc) < 1.38.$ 
To obtain a lower bound for $\a_s$ we choose a  $P$ be in $M_n(\bc)_+$ with $\|P\|_{cbB} =1,$ then, by Lemma \ref{QRnorms},  there exists  a $Q$ in $\cq_n $ such that  Tr$_n(QP) =1.$ By assumption $Q = \a_nR_1 - (\a_n -1) R_2 \leq \a_n R_1$ and we get $$ \|P\|_{cbB} = 1 = \mathrm{Tr}_n(PQ) \leq \a_n\mathrm {Tr}_n(R_1P) \leq \a_n \|P\|_B \leq \a_s\|P\|_B.$$ Hence $k_G^\bc \leq \a_s $ and we have 
\begin{equation}
1.27 < k_G^\bc \leq \a_s \leq  \frac{1}{2-k_G^\bc} < 1.38.
\end{equation}
If the methods from above are applied we get the estimate $K_G^\bc \leq 2 \a_s -1, $ so at its very best we would get $K_G^\bc \leq 2k_G^\bc - 1 < 1.55,$ which is far from Haagerup's upper bound from \cite{Ha3}, which states that $K_G^\bc$ is at most $1.41.$  On the other hand, the inequality above naturally raises the question if $\a_s = k_G^\bc ?$ We think the answer is no, because  Theorem \ref{Chark} gives a characterization of  $k_G^\bc$ which together with the definition of $\a_s$ seems to indicate that $k_G^\bc < \a_s.$

\end{document}